\documentclass[a4paper,10pt]{article}

\usepackage{amssymb,amsfonts,amsmath,geometry}
\usepackage[latin1]{inputenc} 

\newtheorem{theorem}{Theorem}
\newtheorem{corollary}{Corollary}

\newtheorem{proposition}{Proposition}
\newtheorem{definition}{Definition}

\newcommand{\nat}{\mathbb{N}}

\newcommand{\re}{\mathbb{R}}
\newcommand{\cqd}{\, \begin{footnotesize}$\blacksquare$
\end{footnotesize}}

\begin{document}

\title{An Ergodic Theorem on Ergodic Transport}
\author{Joana Mohr and Rafael Rigão Souza} 
\date{\today}

\maketitle

\begin{abstract}
Here we present an ergodic theorem which adapts  a Theorem by J. Elton \cite{E} to the  classical thermodynamical formalism and  to ergodic transport. 
First, we discuss how Elton's theorem can be used to characterise Gibbs measures  for expanding maps. Such characterisation will be done by constructing a stochastic process, defined by a iterated function system (IFS), whose empirical measure converges to the Gibbs measure, in the sense that the mean of any test function evaluated in the outcomes of this stochastic process converges to the integral of such test function with respect to the Gibbs measure. In this way we present a stochastic algorithm that compute integrals of functions.
After this, we turn our attention to ergodic transport: given two sets $X$ and $\Omega$, a measure $\mu$ on $X$ and a dynamics $T$ on $\Omega$, we  consider the set of probability measures on $X \times \Omega$ whose projections on the second coordinate are $T$-invariant, while the projections on the first coordinate are $\mu$. Such measures are called  transport plans. 
 We  call Gibbs plan  any  transport plan that maximizes a pressure functional that is defined  by a potential function added to an entropy term. As in the classical thermodynamical formalism case, we adapt Elton's theorem to define  a stochastic process (using a IFS) whose empirical measures converges to the Gibbs plan.
We provide examples 
and show explicitly calculations in the case where $X$ has two elements and the cost function depends on the two first coordinates of $\Omega$. 

\end{abstract}

\bigskip



\section{Introduction}\label{sec:intro}

Ergodic transport can be motivated by an interesting problem in ergodic theory: given a dynamical system $T:\Omega \to \Omega$, where $(\Omega,d)$ is a compact metric space, $\mathcal{A}$ is the Borel sigma-algebra on $\Omega$  and $T$ is a continuous map, and given a fixed probability measure $\mu$ on $\mathcal{A}$ (which does not need to have any relation to the dynamics $T$), one wants to obtain  the $T$-invariant measure that minimizes the Wasserstein-2 distance to $\mu$. 

If $\nu$ is also a probability measure on $\mathcal{A}$, we denote by $\Pi_{\mu,\nu}$   the set of probability measures on $\Omega \times \Omega$ whose projections on each coordinates are, respectively, $\mu$ and $\nu$. We have that
 the Wasserstein-2 distance between this measures is given by 
$$W_2(\mu,\nu)=\sqrt{\inf_{\pi \in \Pi_{\mu,\nu}}\int d^2(x,y) d\pi(x,y)}.$$
It is known that this distance metrizes the weak-convergence topology on the space of probability measures on $\mathcal{A}$ (see \cite{Vi} for details).

Now, let $\Pi_{\mu,T}$ be the set of probability measures on $\Omega \times \Omega$ whose projections on the first coordinate are equal to $\mu$  and the projections on the second coordinate are $T$-invariant.
If one solves the constrained optimization problem 
$$\inf_{\pi\in \Pi_{\mu,T}} \int d^2(x,y) d\pi(x,y),$$
one obtains (by projecting the solution on the second coordinate) the $T$-invariant measure that is closest to $\mu$ (according to the Wasserstein metric). A measure $\pi$ that solves de minimization problem above is called an  optimal transport plan.

Ergodic transport was studied 
in a more general setting in \cite{LM}. Such paper consider probability measures on $X \times \Omega$,
where $X$ is a set that can be different from  $\Omega$. If $\mu$ is a probability measure on $X$, 
and $T:\Omega \to \Omega$ is a dynamical system defined on $\Omega$,
we denote by $\Pi_{\mu, T}$ the set of probability measures on $X \times \Omega$  that project (in the first coordinate) on  $\mu$, while  projecting, in the second coordinate, on any invariant measure for  $T:\Omega \to 
\Omega$. 
See section \ref{sec:ergodic_transport} for the precise definition of $\Pi_{\mu,T}$. 
 Note that different measures in $\Pi_{\mu, T}$ can project in different invariant measures for $T$.
In \cite{LM}, one is interested in obtaining a measure that attains 
$$\sup_{\pi\in \Pi_{\mu,T}} \int A(x,y) d\pi(x,y),$$
where $A:X \times \Omega \to \re$ is called a potential function. 
If $X$ is a singleton (has only one point), the maximization above is the same found in ergodic optimization problems, and for this reason we can see ergodic transport, in some sense,  as a generalization of ergodic optimization. 
Note that, if we compare the problem above with the one concerning the measure that minimizes Wasserstein measure, we see that we are interested in maximize, rather than minimize, the integral of a function. 
This difference does not pose any conceptual differences in the theory: if one wants to minimize, he just have to change the sign of the potential. 

After  optimal plans were introduced, in \cite{LMMS2} the entropy of plans was defined by means of what can be called the Jacobian of the measure. If we define a functional by adding the integral of the potential with respect to some plan, with the entropy of this plan, and try to maximize this functional, we have what is called a pressure problem (in analogy to classical thermodynamical formalism). One can think that we are considering Thermodynamic Formalism where the potential is random due to the choice of a fixed probability $\mu$.
A variational principle was obtained in \cite{LMMS2} (see section \ref{sec:ergodic_transport} for more details on the results of \cite{LMMS2} that are needed here). 
The plan that satisfies  the variational principle is called an equilibrium plan, and is related to the fixed point of a Ruelle-Perron-Frobenius-type operator (transfer operator).
Moreover, as a relation  between equilibrium plans and optimal  plans, in \cite{LMMS2} it is proved that  optimal plans can be obtained as weak limits of equilibrium plans, if we multiply
the cost by a constant $\beta$, and send $\beta \to \infty$. This is called the zero temperature limit 
(also in analogy to classical thermodynamical formalism - see \cite{Brem}).


In section \ref{sec:ergodic_transport} we review the basic results in ergodic transport, and some details of the results stated above.

Section \ref{sec:Elton_in_ET} brings the main objective of this paper, which is to characterize, via an Ergodic theorem, the equilibrium plan, that solves the pressure problem introduced in  \cite{LMMS2}. This Theorem plays the same role in Ergodic Transport as the Ergodic Theorem in classical Ergodic Theory. 

Such characterisation will be done by constructing a stochastic process  whose empirical measure converges to the Gibbs measure, in the sense that the mean of any test function evaluated in the outcomes of this stochastic process converges to the integral of the test function with respect to the equilibrium plan.
We provide examples and also show explicitly calculations in the case where $X$ has two elements and the cost function depends on the two first coordinates of $\Omega$.

The main tool in getting this caracterization is a result due to Elton, concerning an ergodic theorem for Markov processes defined by iterated function systems (IFS),
which is stated below in an adapted form that is appropriate for our purposes:

\begin{theorem}[Elton - 1987]\label{Elton}
Let Z be a compact metric space, and $\tau_i:Z \to Z$, for $1\leq i \leq d$ be a finite number of Lipschitz contractive maps. 
Let $p_i:Z \to (0,1]$ be Lipschitz continuous weight functions, and suppose $\sum_{i=1}^d p_i(z)=1$ for all $z \in Z$.
Let $\nu$ be a probability measure on  $Z$ that satisfies, for any Borel-measurable set $B$,
\begin{equation}\label{eq}
\nu(B)=\int_Z \sum_{\tau_i(z)\in B}p_i(z) d\nu(z).
\end{equation}

Then, for any  $z_0 \in Z$, if we define by recurrence  
\begin{equation}\label{def-MarkovProcess}
z_{k+1}=\tau_i(z_k) \;\; \mbox{ with probability } \;\; p_i(z_k), 
\end{equation}
we have that, for any continuous function
$f:Z\to \mathbb{R}$, almost surely,
\begin{equation}\label{temp-spac}
	  \frac{1}{N} \sum_{k=0}^{N-1} f(z_k) \to \int f d\nu.
\end{equation}  
\end{theorem}

We will soon explain exactly what almost sure convergence in \eqref{temp-spac} means. Before that, some considerations are necessary: a finite set of maps $\tau_i$ defined on $Z$, together with the {\it transition probabilities} $p_i:Z\to (0,1]$ defines what is called an iterated function system (IFS).  
In this paper we will call  the maps $\tau_i$ as Elton maps.

If we define the transition  Kernel $$P(z,B)=\sum_{\tau_i(z)\in B}p_i(z),$$
\eqref{eq} means that $\nu$ is the invariant measure for the  Markov process  
$\{z_k\}_{k \in \nat}$ defined on $Z$ by \eqref{def-MarkovProcess}.

Elton's theorem implies that the law of large numbers holds for the process $f(z_k)$.

Note that the initial $z_0\in Z$ can be any point in $Z$: the time average on \eqref{temp-spac} does not depend on $z_0$. 


Let us explain the meaning of almost sure convergence in \eqref{temp-spac}:
we begin by denoting $\Omega=\{1,...,d\}^{\nat}$ the Bernoulli set of d symbols, with the $\sigma$-algebra generated by the cylinders.
According to \cite{E}, for any $z \in Z$  there exists (i) a probability $P_z$ on  $\Omega=\{1,...,d\}^{\nat}$, which is given in cylinders by
$$P_z(i_1,i_2,...,i_k)=p_{i_1}(z) p_{i_2}(\tau_{i_1}(z)) 
p_{i_3}(\tau_{i_2}(\tau_{i_1}(z)))... p_{i_k}(\tau_{i_{k-1}}(...(\tau_{i_1}(z))...)), $$  
(ii)  a set $G_z\in \Omega$ such that $P_z(G_z)=1$.
Now, \cite{E} proves that,  if  the {\it address sequence} 
$(i_1,i_2,...)$ 
belongs to $G_z$, then the sequence defined recursively by 
$$\begin{cases} z_0=z, \\
z_{k+1}=\tau_{i_{k+1}}(z_{k}), \mbox{ for } k \geq 0,
\end{cases}$$ is such that \eqref{temp-spac} holds
 for any continuous function $f:Z\to \mathbb{R}$. 
 
In J. Elton´s result the orbits go backward and not forward as in the classical  Birkhoff theorem.



As a simple application of Elton´s result, let $Z=\{1,2,...,d\}$ and $\tau_i=i$. Then $z_k$ defined by \eqref{def-MarkovProcess} is the usual Markov Chain associated to the transition matrix $P_{ij}=p_j(i)$. If $P_{ij}>0$ for all $i,j$, then Elton´s theorem implies, as a particular case, the Law of Large numbers for the Markov chain $z_k$.

Elton's result was proved in a slightly more general form, where the maps $\tau_i$´s only need to be contractive 'on the average', and $Z$ need not to be compact (see details in \cite{E}).  
The invariant measure $\nu$ that satisfies \eqref{eq} was  proved to be unique by some of the authors cited in \cite{E}. More recently a very simple proof of Theorem \ref{Elton} is provided in \cite{FM}, in the case the $p_{i}$´s are constant and $Z$ is compact. 

\medskip


Before  applying Elton´s result to characterize Equilibrium plans in ergodic transport, we will  consider, in section \ref{sec:TF}, an analogous  problem in classical thermodynamical formalism: we know that  the action of the equilibrium measure on test functions can be obtained by a limit procedure considering pre-images of  points (via Ruelle operator - see section \ref{sec:TF}). However, this procedure is not  efficient in computational terms, because the number of preimages growths exponentially  (see section \ref{sec:TF}).   
The method we propose,  using Elton´s theorem, is much more efficient because, using a Markovian stochastic process defined  via an IFS, we get a law of large numbers that gives the integral of test functions as the limit of the mean of such functions  evaluated in the outcomes of the process (see \eqref{temp-spac}).



\section{Ergodic theorem in classical thermodynamical formalism}\label{sec:TF}


In this section we remember the basic facts of classical thermodynamical formalism and then use Elton´s ergodic theorem to characterize Gibbs measures for the shift on the Bernoulli set of symbols.
The results of this section are of independent interest, and will not be used in the following sections.



\bigskip


Let $\Omega=\{1,\hdots,d\}^{\mathbb{N}}$ be the Bernoulli set on $d$ symbols, where $d \in \nat$, $d \geq 2$. 
We know (Tychonoff´s theorem) that $\Omega$ is a compact (metric)  set. We will consider the sigma-algebra generated by the cylinders, which is the Borel sigma-algebra. The dynamics $T$ here is given by the shift map on $\Omega$.

The Ruelle-Perron-Frobenius (see \cite{Ka, OV, PP}) operator (also known as transfer operator) associated to a Lipschitz {\it potential}  $A:\Omega \to \mathbb{R}$ is the operator that associates to any continuous function $\varphi: \Omega \to \re$ the continuous function $L_A(\varphi)$ given by 

$$L_A(\varphi)(z)=\sum_{T(w)=z}e^{A(w)}\varphi(w).$$

A Lipschitz function $A:\Omega \to \re$ is called a {\it normalized potential} if $L_A(1)=1$. It is known (see \cite{PP}) that the RPF operator has a maximal eigenvalue $\lambda_A>0$ associated to an eigenfunction $\varphi_A$, which is Lipschitz, simple and positive (a simple eigenvalue means an eigenvalue that has an associated eigenspace with dimension $1$). 
If $A$ is non-normalized then it can be normalized by considering
$$\bar A = A + \log \varphi_A-\log \varphi_A \circ T - \log \lambda_A.$$

The dual RPF operator, denoted by $L_A^*$, acts on probability measures on $\Omega$, and is defined by $$\int \psi d L_A^*(\mu) = \int L_A(\psi) d\mu.$$

If $\bar A$ is the normalized potential associated to $A$,
we know the dual operator $L_{\bar A}^*$ preserves the simplex of probabilities and therefore has a unique fixed probability $\mu_A$ (i.e. a probability $\mu_A$ that satisfies $L_{\bar A}^*(\mu_A)=\mu_A$), called the Gibbs state associated to $A$, which is invariant and ergodic for $T$, and satisfies 
$$
 \int A d \mu_A + h(\mu_A) = \sup_{ \mu \in \mathcal{M}_T (\Omega)} \left\{ \int A d{\mu} + h(\mu) \right\}.$$
($\mathcal{M}_T (\Omega)$ here means the set of invariant probabilities for $T$.)
The right side of the equation above is called the pressure of $A$. We know that $\mu_A$ is the unique invariant measure to attain the maximum defining the pressure. The unique  measure that attains such  maximum is called the equilibrium measure for $A$, is ergodic, gives positive mass to open sets, and is given by the Gibbs measure $\mu_A$. The last result is also known as the {\it variational principle for pressure}.

Now an interesting question is: how can we characterize the 
Gibbs measure $\mu_A$ ?

More precisely, is it possible to calculate $ \int u d\mu_A$ 
for any test function $u:\Omega\to \re$   ?

A partial answer is given by the fact (see \cite{PP}) that, for any Lipschitz test function $u:\Omega\to \re$, we have
\begin{equation}\label{converg}
    \frac{L_A^n(u)}{\lambda_A^n} \to \varphi_A \int u d\mu_A\,,
\end{equation}
when $n \to \infty$, where 
$$L_A^n(\varphi)(z)=\sum_{T^n(w)=z}e^{S_n(A)(w)}\varphi(w),$$
 where $S_n(A)=\sum_{k=0}^{n-1} A \circ T^k$ is the Birkhoff sum  of order $n$ of $A$. 
The convergence above is on the uniform convergence topology.  

However, this limit procedure above is not of practical use because it involves the evaluation of $u$ in the inverse images of order $n$, with $n \to \infty$, of some point of $\Omega$, and such calculation is not practical to be implemented: for the shift on $d$ symbols, we have $d^n$ inverse images of order $n$. To get things even worse, we have to evaluate the Birkhoff sum of the potential $A$ of order $n$ in each one of the $d^n$ inverse images of a chosen point.

Therefore, it is necessary another way for calculating $ \int u d\mu_A$.
This can be accomplished  by using Elton´s theorem. In this way we get a much  more efficient procedure, because what we will get is a Markovian stochastic process (that can be easily simulated in a Monte-Carlo process) and we will get a  law of large numbers that will result in the integral of any test function, as shown in equation \eqref{temp-spac}.

Suppose $A:\Omega\to \re$ is normalized, i. e. $\sum_{i=1}^d e^{A(iz)}=1 \;\; \forall z \in \Omega$.

Let  $\tau_i(z)=iz$ be the Elton maps, and the transition probabilities be  given by $p_i(z)=e^{A(iz)}$. The normalization hypothesis on $A$ implies that $\sum p_i(z)=1$ for any $z$, and also the Lipschitz continuity of
$A$ implies that $p_i$ is Lipschitz. 


Let $\mu$ be the probability measure that satisfies  \eqref{eq} in Elton´s theorem.
Let $B$ be any measurable set of  $\Omega$:

Then $$\mu(B)= 
\int \sum_{\tau_i(z)\in B}p_i(z) d\mu(z) = 
\int \sum_{iz \in B} e^{A(iz)} d\mu(z) =$$
$$ = \int \sum_{i} e^{A(iz)} I_B(iz) d\mu(z)
= \int \mathcal{L}_A(I_B) d\mu= \mathcal{L}_A^*\mu(B),$$
and therefore $ \mu$ is the fixed point of Ruelle dual operator.

As a result, Elton measure $\mu$ coincides with the Gibbs measure $\mu_A$, and we have 
 
\begin{theorem}[Birkhoff-Elton Theorem in classical thermodynamical formalism]

If we choose any  $z_0 \in \Omega$, and then, by recurrence, choose 
$$z_{k+1}=iz_k \; \mbox{ with probability } p_i(z_k)=e^{A(iz_k)},$$ then, for any continuous test function
$f:\Omega\to \mathbb{R}$,
$$
\frac{1}{N} \sum_{k=0}^{N-1} f(z_k) \to \int f d \mu_A,
$$
with probability one. 

\end{theorem}

\section{Ergodic Transport}\label{sec:ergodic_transport}

Now we recall the main concepts of Ergodic Transport. See \cite{LMMS2} for more details. 

Let $X$ be a finite set and $\Omega=\{1,...,d\}^{\mathbb{N}}$ the Bernoulli space on $d$ symbols.

\bigskip

We denote the RPF operator associated to a Lipschitz cost (potential) $c:X\times \Omega\to \mathbb{R }$ as 
\[
 L_c(v)(y) = \sum_{x\in X}\sum_{\sigma(w)=y} e^{c(x,w)}v(w)=\sum_{\sigma(w)=y}\Big(\sum_{x\in X} e^{c(x,w)}\Big)v(w),
 \]
 for any  $v\in C(\Omega)$. We remark that this is the classical RPF operator associated to the potential $b_c(y)=\log (\sum_{x\in X} e^{c(x,y)})$. We know that $L_c$ has a maximal eigenvalue $\lambda_c$, which is simple and positive, and there is a positive eigenfunction $h_c$  associated to $\lambda_c$, see \cite{PP}.
 
  The  RPF extended to the continuous  functions   $X \times  \Omega$ is defined as 
\[
\hat L_c(u)(y) = \sum_{x\in X}\sum_{\sigma(w)=y} e^{c(x,w)}u(x,w),
\] 
for any $u:X\times \Omega \to \mathbb{R}$. 

Note that $\hat L_c$ sends $u:X \times  \Omega\to \mathbb{R}$ to the function denoted by $\hat L_c(u):\Omega \to \mathbb{R}$.

\begin{definition}\label{normalization} We say that a Lipschitz cost (potential) $c:X\times \Omega\to \mathbb{R }$ is normalized if for any $y\in \Omega$, we have
\[\sum_{x\in X} \sum_{\sigma(w)=y}e^{c(x,w)} = 1.\]
\end{definition}

If $c$ is a Lipschitz cost that is not normalized, we associate to $c$ the normalized cost 
\begin{equation}\label{cnormal}
\bar{c}(x,y)=c(x,y)+\log(h_c(y))-\log(h_c\circ \sigma(y))-\log(\lambda_c),
\end{equation} were $h_c$ is the positive eigenfunction associated to the maximal eigenvalue $\lambda_c$ of $L_c$.

\medskip

If $c$ is  normalized,  we denote  $\hat L_c^*$
 the operator on $P(X\times \Omega)$ defined by 
\begin{equation}\label{Op}
 \hat L_c^*(\pi)(u(x,y))\,=\,\int_{X\times \Omega} \left(\sum_{\alpha\in X}\sum_{\sigma(w)=y} e^{c(\alpha,w)} u(\alpha,w) \right) \, d\pi(x,y)\,.
\end{equation}


The normalization property implies that $\hat L_c^*(\pi)(1)=1$, i.e.
$\hat L_c^*$ preserves the convex and compact set $P(X\times \Omega)$, and Tychonoff-Schauder theorem implies

\begin{proposition} \label{Gi} Given a  normalized cost $c$ there exists a {unique} fixed point for the operator $\hat L_c^*$. It will be denoted by $\pi_c$.
\end{proposition}

Therefore, we have 
 \begin{equation}\label{pic}
 \pi_c(u(x,y))\,=\,\int_{X\times \Omega} \left(\sum_{\alpha\in X}\sum_{\sigma(w)=y} e^{c(\alpha,w)} u(\alpha,w) \right) \, d\pi_c(x,y).
 \end{equation}

\begin{definition}
The fixed point, $\pi_c$, for  $\hat L_c^*$ is called the {\bf  Gibbs plan} for the normalized cost (potential) $c$.

\end{definition}

We denote by $\Pi(\cdot,\sigma)$ the set of plans such that its $y$-marginal is $\sigma$-invariant, i.e., 
\begin{equation}\label{pipontosigma}
\int_{X\times \Omega}  g(y) \,d\pi(x,y) = \int_{X\times \Omega}  g(\sigma(y)) \, d\pi(x,y) \ \ \text{for any} \, g\in C(\Omega).
\end{equation}

By Theorem 4 of \cite{LMMS2}, we have that  $\pi_c\in\Pi(\cdot,\sigma)$  and the $y$-marginal of $\pi_c$ is the  Gibbs measure $\nu_c$  to the classical RPF operator with the potential $b_c(y)=\log (\sum_{x\in X} e^{c(x,y)})$.

 \medskip



We denote by $[x,y_1...y_n]= \{(\alpha,w)\in X\times\Omega:\alpha=x,w_1=y_1,...,w_n=y_n\}$ and $[y_2...y_n]=\{w\in \Omega:w_1=y_2,...,w_{n-1}=y_n\}. $
Consider a fixed  plan  $\pi\in \Pi(\cdot,\sigma)$ with $y$-marginal $\nu$ and define $$J_\pi^{n}(x,y) = \frac{\pi([x,y_1...y_n])}{\nu([y_2...y_n])}$$ if $y=(y_2,y_3,...) \in \text{ supp}(\nu)$. From the Increasing Martingale Theorem the functions $J_\pi^{n}$ converge to a function $J_\pi(x,y)$ in $L^{1}(X\times\Omega,\mathcal{A}(X\times\Omega),\pi)$ and for $\pi$ a.e. $(x,y)$. For each plan $\pi$ this function $J_\pi$ can be also obtained via the Radon-Nikodyn Theorem. 
 We have, $J_\pi>0$ a.e. ($\pi$) and
$ \sum_{x\in X}\sum_{a\in \{1,..,d\}} J_\pi(x,ay) = 1$.

We define the entropy of a plan $\pi$ as 
$$H(\pi) = -\int \log (J_{\pi})\, d\pi .$$

\begin{definition}\label{def:pressure}
The pressure of a Lipschitz continuous  cost (potential) $c$ is defined by
\[P(c) = \sup_{\pi \in \Pi(\cdot,\sigma)} \left(\int_{X\times\Omega} c\, d\pi + H(\pi)\right).\]

A plan $\pi\in \Pi(\cdot,\sigma) $ which realizes the supremum is called {\bf an equilibrium plan for $c$}.
\end{definition}

\begin{theorem}[Variational Principle over $\Pi(\cdot,\sigma)$]\label{gibbs-equilibrium}

 Let us fix a Lipschitz cost $c$. Then,
$P(c) =\log(\lambda_c)$, where $\lambda_c$ is the main eigenvalue of $L_c$. The equilibrium plan for $c$ is { unique} and given by the Gibbs plan for the normalized cost $\overline{c}:=c +\log(h_c) - \log(h_c\circ\sigma) - \log(\lambda_c)$, where $h_c$ is the positive eigenfunction associated to $\lambda_c$.

\end{theorem}

\bigskip

Now let us fix $\mu$ a probability on $X$. We define the $\mu$-pressure of $c$ by

\begin{equation}\label{mu-pressure}
 P_{\mu}(c)=\sup_{\pi \in \Pi(\mu,\sigma)}  \int_{X\times \Omega} c \, d\pi + H(\pi),
\end{equation}
where $\Pi(\mu,\sigma)$ is the set of all plans satisfying
\begin{equation}\label{pimusigma}
\left\{ \begin{array}{l}
\int_{X\times \Omega} f(x) \, d\pi(x,y) = \int_X f(x) \, d\mu(x) \ \ \text{for any} \, f\in C(X), \\
\int_{X\times \Omega}  g(y) \,d\pi(x,y) = \int_{X\times \Omega}  g(\sigma(y)) \, d\pi(x,y) \ \ \text{for any} \, g\in C(\Omega),
\end{array}
\right.
\end{equation}
which means, the set of probabilities $\pi$ such that the $x$-marginal of $\pi$ is the fixed probability $\mu \in P(X)$ and the $y$-marginal of $\pi$ is $\sigma$-invariant.

Note that $ P_{\mu}(c)\leq  P(c)$.

By compactness, there exists a  plan $\tilde\pi_c\in\Pi(\mu,\sigma)$ which attains the supremum at \eqref{mu-pressure}.

We have the following duality result:

\begin{theorem}\label{dualidade} 

Given a Lipschitz cost $c$, and let $\mu$ be a probability on $X$, we have
\begin{equation}\label{Fenchel-Rockafellar equation}
\,\,\,\,\
 \,\inf_{\varphi:P(c-\varphi)= 0}  \int_X \varphi \, d\mu =  \sup_{\pi \in \Pi(\mu,\sigma)}\int_{X\times \Omega}\, c(x,y)\, d \pi + H(\pi).
\end{equation}
The minimization is performed on the set of functions $\varphi:X \to \re$ such that $P(c-\varphi)=0$. 
The supremum at ($\ref{Fenchel-Rockafellar equation}$) is attained in at least one plan, while the infimum is attained in exactly one function $\tilde \varphi$.

\end{theorem}
\medskip

Moreover, if $\tilde{\varphi}:X \to \re$ is the unique minimizer of \eqref{Fenchel-Rockafellar equation}, then $P(c-\tilde{\varphi})=0$ and also
$$  \int_X \tilde{\varphi} \, d\mu =  \sup_{\pi \in \Pi(\mu,\sigma)}\int_{X\times \Omega}\, c(x,y)\, d \pi + H(\pi)$$
which implies
$$  \sup_{\pi \in \Pi(\mu,\sigma)}\int_{X\times \Omega}\, (c(x,y)-\tilde{\varphi} )\, d \pi + H(\pi)=0,$$
which means $P_{\mu}(c-\tilde \varphi)=0$.

\medskip

Therefore $P(c-\tilde \varphi)=P_{\mu}(c-\tilde \varphi)=0$, and the maximizer of \eqref{Fenchel-Rockafellar equation} is the equilibrium plan for $c-\tilde \varphi$, i.e. satisfies the non-constrained variational principle. We just proved the

\begin{corollary}\label{corolmin}
Let $\tilde{\varphi}:X \to \re$ be the unique minimizer for the
Fenchel-Rockafellar duality
\eqref{Fenchel-Rockafellar equation}. The equilibrium plan $\pi_{c- \tilde{\varphi}}$,
for $c- \tilde{\varphi}$, belongs to $\Pi(\mu,\sigma)$ and is the unique maximizer of \eqref{Fenchel-Rockafellar equation}.
\end{corollary}

As a final remark in this section, let us observe that  we show, in \cite{LMMS2}, among other things, how equilibrium plans can be used to obtain the optimal transport plan, by a limit procedure where the cost $c$ is multiplied by a constant $\beta$ (the inverse of the temperature) and the parameter $\beta\to \infty$.


\vspace{1cm}


\section{An ergodic theorem in ergodic transport}\label{sec:Elton_in_ET}

Now we discuss the main goal of this paper: the characterization, via Elton theorem, of the plan $\pi_{c,\mu}$ that satisfies the supremum
$$\sup_{\pi \in \Pi(\mu,\sigma)}\int_{X\times \Omega}\, c(x,y)\, d \pi + H(\pi).
$$
where $X=\{1,2\}$, $\mu$ is a probability measure on $X$, $c:X \times \Omega \to \mathbb{R}$ depends  on $x\in X$ and only on the two first coordinates of $y\in\Omega$, 
and $\Pi(\mu,\sigma)$ is given by \eqref{pimusigma}. Note that we are dealing with the {\bf constrained} optimization problem, i.e. we search an optimising measure among those measures that project on a fixed probability measure $\mu$ on $X$. Any such measure is given by $\mu=(p, 1-p)$, where $0\leq p \leq 1$. 

\vspace{0.3cm}

By the corollary \ref{corolmin}, the plan $\pi_{c,\mu}$  is the equilibrium plan $\pi_{c- \tilde{\varphi}}$, where  $\tilde{\varphi}$ is the minimizer of  \eqref{Fenchel-Rockafellar equation}.
Note that this minimization is performed on the set of functions $\varphi:X \to \re$ such that $P(c-\varphi)=0$. 

\bigskip
We will do such characterisation in three steps:

\bigskip {\bf Step 1} We need to find the minimizer $\tilde \varphi = (\varphi_1,\varphi_2) \in \mathbb{R}^2$ that solves the minimization problem \eqref{Fenchel-Rockafellar equation}, i.e, that solves
$$\inf_{\varphi:P(c-\varphi)= 0}  \int_X \varphi \, d\mu.$$ 
Such minimizer exists and is unique (theorem \ref{dualidade}).

\bigskip {\bf Step 2} If $A= c - \tilde{\varphi}$ , we need to find the maximal eigenvalue $\lambda_A$  for the RPF operator $ L_A$, as well as the associated eigenfunction $h_A$. Once we have that, we define the normalized potential 
$\bar A= A +\log(h_A)-\log(h_A \circ \sigma) - \log(\lambda_A),$
which can be reduced to 
$$\bar A= A +\log(h_A)-\log(h_A \circ \sigma),$$
as $\log(\lambda_A)=P(A)=P(c-\tilde \varphi)=0$.

\bigskip {\bf Step 3}
We use Elton Theorem to characterize the constrained equilibrium measure $\pi_{c,\mu}$ showing how one can get the integral of any test function by means of a stochastic simulation (a Monte Carlo method).
More precisely, for any continuous function  $f:X\times \Omega \to \mathbb{R}$, show that, almost certainly,  $$\frac{1}{N} \sum_{k=0}^{N-1} f(x^k,y^k) \to \int f d\pi_{c,\mu},$$
for a well-choosen sequence $(x^k,y^k)$, obtained by a probabilistic algorithm.


\bigskip

\bigskip

\textbf{Analysis of Step 3:}

Now we suppose Steps 1 and 2 were successfully performed (the analysis of these two first steps will be done after theorem \ref{Birkoff-Elton}).
We proceed to step 3: 

\medskip

By Corollary \ref{corolmin} we know that $\pi_{c,\mu}$ is the equilibrium measure $\pi_A$ or $\pi_{c- \tilde{\varphi}}$, and by theorem \ref{gibbs-equilibrium} we have that $\pi_A$ or $\pi_{c- \tilde{\varphi}}$ is the Gibbs plan for $\bar A$, i.e., $$\hat{L}^*_{\bar A}(\pi_A)=\pi_A.$$

Now we want to show that $\pi_A$ (i.e. $\pi_{c,\mu}$) satisfies the hypothesis of theorem \ref{Elton}.


Let us take, following the notation used in Elton's Theorem \ref{Elton},  $Z=X \times \Omega$,  then we define for all
 $(\alpha,i)\in X\times\{1,...,d\}$ and for all
 $(x,y) \in Z$, the weight function
 $$p_{\alpha,i}(x,y)=e^{\bar A(\alpha ,i y)}.$$ Note that $p_{\alpha,i}$ depends only on $y$.

As $\bar A$ is normalized we have
$$1=\sum_{\alpha \in X}\sum_{ 1 \leq i \leq d}e^{\bar A(\alpha ,i y)}=\sum_{\alpha \in X, 1 \leq i \leq d} p_{\alpha,i}(x,y)\,\,, \;\; \forall  (x,y) \in Z.$$ 

For each $\alpha \in X, 1 \leq i \leq d$, we define the  Elton map
$$\tau_{\alpha,i}(x,y)=(\alpha ,i y).
$$

Finally, if  $B$ is a borelian set of $Z = X\times \Omega$, and $\chi_B$ is the  characteristic function of $B$. We have, by equation \eqref{pic} that

$$\pi_A(B)=\int \hat{L}_{\bar A} (\chi_B)d \pi_A=
\int \sum_{\alpha \in X, 1 \leq i \leq d} p_{\alpha,i}(x,y) \chi_B(\alpha ,i y) d \pi_A(x,y)
$$
$$=
\int \sum_{\tau_{\alpha,i}(x,y) \in B} p_{\alpha,i}(x,y)  d \pi_A(x,y),$$
which  means that 
 $\pi_A$ satisfies  \eqref{eq} and the Elton Theorem is true for the stochastic process generated by   $p_{\alpha,i}$ and
$\tau_{\alpha,i}$. Therefore, we have proved the following theorem:

\begin{theorem}\label{Birkoff-Elton}[Birkhoff-Elton Theorem in Ergodic Transport]
Fix any $(x^0,y^0)\in X\times \Omega$. Define by recurrence  $$(x^{k+1},y^{k+1})=(\alpha ,i y^k) \;\mbox{ with probability }
 \;e^{\bar A(\alpha ,i y^k)},$$ where $\bar A$ is the normalized cost associated to $A=c - \tilde\varphi$, and
 $\tilde{\varphi}$ is the only minimizer of the Fenchel-Rockafellar duality equation  \eqref{Fenchel-Rockafellar equation}.
Then, for any continuous function  $f:X\times \Omega \to \mathbb{R}$, we have, almost certainly,  $$\frac{1}{N} \sum_{k=0}^{N-1} f(x^k,y^k) \to \int f d\pi_A,$$
where  $\pi_A=\pi_{c,\mu}$ is the constrained maximizer of the duality equation  \eqref{Fenchel-Rockafellar equation}.

\end{theorem}

\textbf{Analysis of Steps 1 and 2:}

We will consider here the case $d=2$, i.e., $\Omega=\{1,2\}^{\mathbb{N}}$. As  $c$   depends on $x$ and only on the two first coordinates of $y$, that is, $c(x,y) = c(x,\,y_1,\,y_2)$, it is represented by the following matrices:

\begin{equation}\label{defmatrix}
C^{1}=\left(\begin{array}{cc}e^{c^1_{11}} & e^{c^1_{12}}\\ e^{c^1_{21}} & e^{c^1_{22}} \end{array}\right)\,\,\,\,,\,\,\,\,C^{2}=\left(\begin{array}{cc}e^{c^2_{11}} & e^{c^2_{12}}\\ e^{c^2_{21}} & e^{c^2_{22}} \end{array}\right),
\end{equation}

where $C_{ij}^x=e^{c^x_{ij}}= e^{ c(x,i,j)}$, $x,i,j=1,2$.

Note that the operator $L_c:C(\Omega)\to C(\Omega)$, for any  $v: \Omega \to \mathbb{R}$, can be rewritten as
$$L_c(v)(y)=\sum_{i\in \{1,2\}}\Big(\sum_{x\in X} e^{c(x,iy)}\Big)v(iy)=\sum_{i\in \{1,2\}}e^{b_c(iy)}v(iy),$$
 where $b_c(y)=\log \Big(\sum_{x\in X} e^{c(x,y)}\Big)$. Then $b_c:\Omega\to \mathbb{R}$ depends on the two first coordinates of $y$ and is represented by the matrix 
\begin{equation}\label{defB}
B=\left(\begin{array}{cc}         e^{c^{1}_{11}} \, +e^{c^{2}_{11}} & e^{c^{1}_{12}} \, +e^{c^{2}_{12}}\,\\ e^{c^{1}_{21}}\,+ e^{c^{2}_{21}}\,& e^{c^{1}_{22}}\,+ e^{c^{2}_{22}}\, \end{array}\right).
\end{equation}
It is a well known fact, see \cite{PP}, that the positive eigenfunction, that we denote by $v_c$, associated to the main eigenvalue, $\lambda_c$,  of $L_c$ depends only on the first coordinate of $y$, and this means that $v_c$ is in fact a vector of $\mathbb{R}^2$. 

More precisely, $v_c$ is the  left eigenvector of $B$ associated to the eigenvalue $\lambda_c$, because
\begin{equation}\label{autofuncao}
L_c(v_c)(y)=\lambda_c v_c(y)\iff v_cB=\lambda_c v_c.
\end{equation}
 
\bigskip

Let also  
\[C^{12}=\left(\begin{array}{cc}e^{c^1_{11}} & e^{c^1_{12}}\\e^{c^2_{21}} & e^{c^2_{22}} \end{array}\right)\,\,\,\,,\,\,\,C^{21}=\left(\begin{array}{cc}e^{c^2_{11}} & e^{c^2_{12}}\\ e^{c^1_{21}} & e^{c^1_{22}} \end{array}\right).\]

\bigskip
The next Theorem shows how to perform the Step 1.

\begin{theorem}\label{varphi-custoduascoordenadas}
Let $\mu=(p,1-p)$. Suppose $c:X\times \Omega\to \mathbb{R}$ is a Lipschitz cost that depends on $x$ and just on the two first coordinates of $y$. 

(a) 
Then the unique minimizer $\tilde{\varphi}$ of the left-hand side of  the equation
\begin{equation}\label{eqdualidade}
\inf_{\varphi:P(c-\varphi)= 0} \int_{X} \varphi(x) \, d\mu =
 \sup_{\pi \in \Pi(\mu,\sigma)} \int_{X\times \Omega}  c  \, d\pi\, +\, H(\pi)
\end{equation}
 is given by  
$$\begin{cases}
\tilde\varphi_1=-\log (z_1), \\ 
\tilde\varphi_2=-\log (z_2),
\end{cases}$$
where $(z_1,z_2)$ belongs to the set of positive solutions of the following system of equations
\begin{equation}\label{eqML}
\begin{cases}
a\,z_1^2+b\,z_2^2+c\,z_1z_2+d\,z_1+e\,z_2+1=0, \\ 
 2bpz_2^2-2a(1-p)z_1^2+c(2p-1)z_1z_2+epz_2-d(1-p)z_1= 0 ,
\end{cases}
\end{equation}
where $a=\det C^1$, $b=\det C^2$, 
$c=\det C^{12}+\det C^{21}$, 
$d=-\mbox{tr}\, C^1$,  and $e=-\mbox{tr}\, C^2$.

This means that $(z_1,z_2)$ belongs to the set of intersection points of the two conics given by the equations above. Such intersections points are four, at most, and 
$(z_1,z_2)$ is the one that
 minimizes $-p \log (z_1) - (1-p) \log (z_2)$.




(b) Suppose also that $C^1$ and $C^2$ are stochastic matrices. 
Then the unique minimizer $\tilde{\varphi}$ of the left-hand side of  the equation \eqref{eqdualidade}
is given by 
$$\begin{cases}
\tilde{\varphi}_1=-\log (p), \\ 
\tilde{\varphi}_2=-\log (1-p).
\end{cases}$$
Also we have
$$ \sup_{\pi \in \Pi(\mu,\sigma)} \int_{X\times \Omega}  c  \, d\pi\, +\, H(\pi)= -p \log (p) - (1-p) \log (1-p).$$
\end{theorem}

{\bf Proof:}
First we need to understand what means the restriction, $P(c-\varphi)= 0$, in left-hand side of equation \eqref{eqdualidade}.
By Theorem \ref{gibbs-equilibrium} $P(c-\varphi)=\log(\lambda_{c-\varphi})$, where $\lambda_{c-\varphi}$ is the main eigenvalue of $L_{c-\varphi}$. Therefore $P(c-\varphi)=0$ means $\lambda_{c-\varphi}=1$.
As $\varphi:X\to\mathbb{R}$ and $X=\{1,2\}$, $\varphi$ can be seen as two-dimensional vector $\varphi=(\varphi_1,\varphi_2)$. Then we have

$$\displaystyle L_{c-\varphi}\psi(y)=\sum_{i}\Big[e^{c(1,ij)-\varphi(1)}+e^{c(2,ij)-\varphi(2)}\Big]\psi(iy),$$ where $y=(j,y_2,y_3,...)$.
As in \eqref{defB}, the operator $L_{c-\varphi}$ is represented by the $2 \times 2$ matrix 

\begin{equation}\label{dominant}
B=\left(\begin{array}{cc}         e^{c^{1}_{11}} \,e^{-\varphi_1} +e^{c^{2}_{11}}\,e^{-\varphi_2} & e^{c^{1}_{12}} \,e^{-\varphi_1} +e^{c^{2}_{12}}\,e^{-\varphi_2}\\ e^{c^{1}_{21}}\,e^{-\varphi_1}+ e^{c^{2}_{21}}\,e^{-\varphi_2} & e^{c^{1}_{22}}\,e^{-\varphi_1}+ e^{c^{2}_{22}}\,e^{-\varphi_2} \end{array}\right)
\end{equation}
finally using \eqref{autofuncao}, $P(c-\varphi)=0$ means   $\lambda_{c-\varphi}=1$ is the dominant   eigenvalue of $B$.

If we apply the change of coordinates  $z_1= e^{-\varphi_1}$, $z_2= e^{-\varphi_2}$.  
Then  $P(c-\varphi)=0$ if and only if $(z_1,z_2)$ satisfy 
\begin{enumerate}
\item[(i)] $z_1>0,z_2>0$,
\item[(ii)] 
$$\det \, \left(\begin{array}{cc}         e^{c^{1}_{11}} \,z_1   +e^{c^{2}_{11}}\,z_2\,-1 & e^{c^{1}_{12}} \,z_1 +e^{c^{2}_{12}}\,z_2\\ e^{c^{1}_{21}}\,z_1+ e^{c^{2}_{21}}\,z_2 & e^{c^{1}_{22}}\,z_1+ e^{c^{2}_{22}}\,z_2-1 \end{array}\right)\,=0,$$
\item[(iii)]  the other eigenvalue of the matrix $B$ given in \eqref{dominant} is less than 1.
\end{enumerate}

Note that we want to find the unique minimizer of 
$$\inf_{\varphi:P(c-\varphi)= 0} \int_{X} \varphi(x) \, d\mu ,$$
but instead of that we will minimize $\int_{X} \varphi(x) \, d\mu ,$ subject to the restriction that $(z_1,z_2)$  satisfies condition (ii)\footnote{see remark after Theorem \ref{varphi-custoduascoordenadas}}. After we get all the solutions of this problem we will test which of them also satisfies conditions (i) and (iii).

\medskip


 Condition (ii) is equivalent to 
$$(e^{c^{1}_{11}} \,z_1   +e^{c^{2}_{11}}\,z_2\,-1)(e^{c^{1}_{22}}\,z_1+ e^{c^{2}_{22}}\,z_2-1)-(e^{c^{1}_{12}} \,z_1 +e^{c^{2}_{12}}\,z_2)( e^{c^{1}_{21}}\,z_1+ e^{c^{2}_{21}}\,z_2)  =$$
$$= z_1^2(e^{c^{1}_{11}}e^{c^{1}_{22}}-e^{c^{1}_{12}}e^{c^{1}_{21}})+z_2^2(e^{c^{2}_{11}}e^{c^{2}_{22}}-e^{c^{2}_{12}}e^{c^{2}_{21}})+z_1(-e^{c^{1}_{11}}-e^{c^{1}_{22}})+z_2(-e^{c^{2}_{11}}-e^{c^{2}_{22}})+$$
$$+z_1z_2(e^{c^{1}_{11}}e^{c^{2}_{22}}+e^{c^{2}_{11}}e^{c^{1}_{22}}-e^{c^{1}_{12}}e^{c^{2}_{21}}-e^{c^{2}_{12}}e^{c^{1}_{21}})+1= $$
$$:=a\,z_1^2+b\,z_2^2+c\,z_1z_2+d\,z_1+e\,z_2+1=0.$$




If we denote by
$$g(z_1,z_2)=a\,z_1^2+b\,z_2^2+c\,z_1z_2+d\,z_1+e\,z_2+1$$ 
then  condition (ii) describes an algebraic curve (a conic) on $\mathbb{R}^2$, given  by $$g(z_1,z_2)=0.$$ 


Therefore we need to  minimize  $$f(z_1,z_2):=\int_X\varphi(x) d\mu=\varphi_1p+\varphi_2 (1-p)=-p\log z_1 -(1-p)\log z_2,$$ subject to the restriction  $g(z_1,z_2)=0$.
To do that we will use the Lagrange multiplier theorem: 
if $(z_1,z_2)$ is a solution to the constrained minimization problem above, there exists a $\lambda\in \mathbb{R}$ such that $$\nabla f(z_1,z_2)=\lambda \nabla g(z_1,z_2).$$

Using that $$ \begin{cases} \nabla g(z_1,z_2)=(2az_1+cz_2+d,2bz_2+cz_1+e), \\
\nabla f(z_1,z_2)=(-\frac{p}{z_1},-\frac{1-p}{z_2}), \end{cases}$$
we need to solve the following equations 
$$-\frac{p}{z_1}=\lambda(2az_1+cz_2+d) \,\,\,\mbox{and}\,\,\, -\frac{1-p}{z_2}=\lambda (2bz_2+cz_1+e), $$
or
$$\frac{p}{z_1(2az_1+cz_2+d)}= \frac{1-p}{z_2(2bz_2+cz_1+e)}, $$
or
$${p (2bz_2^2+cz_1z_2+ez_2)}= {(1-p) (2az_1^2+cz_2z_1+dz_1)}, $$
which are equivalent to
\begin{equation}\label{eqml0}
 2bpz_2^2-2a(1-p)z_1^2+c(2p-1)z_1z_2+epz_2-d(1-p)z_1= 0 
\end{equation}
which, together with equation $g=0$ give the system of equations \eqref{eqML}.
 
  In order to conclude the proof of item (a) of Theorem \ref{varphi-custoduascoordenadas}, we remark that the unique minimizer of $\displaystyle\inf_{\varphi:P(c-\varphi)= 0} \int_{X} \varphi(x) \, d\mu$ belongs to the intersection of the two conics in \eqref{eqML}, which are, at most, four points. For each of this solutions we test if it satisfies conditions (i), (that is, the solution have to be in the positive quadrant) and (iii). If there are more than one solution  that satisfies conditions (i)-(iii), then we test which one is the minimizer.  






 \medskip

Now, in order to prove item (b) of Theorem \ref{varphi-custoduascoordenadas},
 let us suppose that $C^1$ and $C^2$ are stochastic matrices, i.e.
\begin{equation}\label{stochmatr}
C^{1}=\left(\begin{array}{cc}e^{c^1_{11}} & e^{c^1_{12}}\\ 1-e^{c^1_{11}} & 1- e^{c^1_{12}} \end{array}\right)\,\,\,,\,\,\,C^{2}=\left(\begin{array}{cc}e^{c^2_{11}} & e^{c^2_{12}}\\ 1-e^{c^2_{11}} & 1- e^{c^2_{12}} \end{array}\right)
\end{equation}
We have $a=e^{c^1_{11}}-e^{c^1_{12}}$, $b=e^{c^2_{11}}-e^{c^2_{12}}$, $c=a+b$, $d=-(a+1)$, $e=-(b+1)$,

and equation $g=0$ becomes 
\begin{equation}\label{g0}
a\,z_1^2+b\,z_2^2+(a+b)\,z_1z_2-(a+1)\,z_1-(b+1)\,z_2+1=0
\end{equation} while equation \eqref{eqml0} becomes
\begin{equation}\label{eqml}
 2bpz_2^2-2a(1-p)z_1^2+(a+b)(2p-1)z_1z_2-(b+1)pz_2+(a+1)(1-p)z_1= 0 
\end{equation}

First we solve equation \eqref{g0} in terms of $z_2$. We get, if $a\neq 0$,  two solutions 
$\tilde z_1=\frac{1-bz_2}{a}$ and $\bar z_1=1-z_2$, and if $a=0$ we get only the solution $\bar z_1=1-z_2$.

Now, we solve \eqref{eqml}: 

  
Case 1: if $\tilde z_1=\frac{1-bz_2}{a}$, \eqref{eqml} becomes 

$$  2bpz_2^2-2a(1-p)\bigg(\frac{1-bz_2}{a}\bigg)^2+(a+b)(2p-1)\frac{1-bz_2}{a}z_2-(b+1)pz_2+$$
$$+(a+1)(1-p)\frac{1-bz_2}{a}= 0.$$
This equation can be rewritten as
$$\bigg(b-\frac{b^2}{a}   \bigg) z^2_2+\bigg( \frac{2b-bp}{a}+p-1-b \bigg)z_2 +\frac{p-1}{a}+(1-p) = 0$$
or 
$$(ba-b^2   ) z^2_2+( 2b-bp+pa-a-ba )z_2 +p-1+a-pa = 0$$
finally, this is equivalent to 

$$\Big(a(z_2-1)-bz_2+1\Big)\Big(bz_2+p-1\Big)=0, $$
which  have two solutions 
$\tilde z_2=\frac{a-1}{a-b}$ and $\bar z_2=\frac{1-p}{b}$. 

Then, we get two possible pair of  solutions $(z_1,z_2)$ for the system given by \eqref{g0} and \eqref{eqml}: if $\tilde z_2=\frac{a-1}{a-b}$, then $\tilde z_1=\frac{1-bz_2}{a}=\frac{1-b}{a-b}$
and if 
$\bar z_2=\frac{1-p}{b}$, then $\tilde z_1=\frac{1-bz_2}{a}=\frac{p}{a}$.

Case 2:  $\bar{z_1}=1-z_2$, then equation \eqref{eqml} becomes

$$2bpz_2^2-2a(1-p)(1-2z_2+z_2^2)+(a+b)(2p-1)(z_2-z_2^2)-(b+1)pz_2+$$
$$+(a+1)(1-p)(1-z_2)= 0$$



which can be rewritten as

$$\Big(-a+b \Big)z_2^2+\Big(2a-ap+bp-b-1\Big)z_2+\Big(ap-a-p+1\Big)=0.$$

Finally, this is equivalent to 

$$\Big(p+z_2-1 \Big)  \Big(a(z_2-1)-bz_2+1 \Big)=0$$
which  have two solutions 
$\tilde z_2=\frac{a-1}{a-b}$ and $\bar z_2=1-p$. 

Then, we get another two possible pair of  solutions $(z_1,z_2)$:

If $\tilde z_2=\frac{a-1}{a-b}$, then $\tilde z_1=1-z_2=\frac{1-b}{a-b}$
and if $\bar z_2=1-p$, then $\tilde z_1=1-z_2=p$.

Collecting all the  solutions of cases 1 and 2, we have 3 different solutions $(z_1,z_2)$ for the system    given by \eqref{g0} and \eqref{eqml} 
$$\bigg(\frac{1-b}{a-b},\frac{a-1}{a-b}\bigg)\;\;\;, \;\;\;\bigg(\frac{p}{a},\frac{1-p}{b}\bigg)\;\;\;, \;\;\; (p,1-p).$$ 

Now we test conditions (i) and (iii) for each one of the three possible solutions above:


\medskip

\textbf{Claim:} $\big(\frac{1-b}{a-b},\frac{a-1}{a-b}\big)$ do not satisfy item (i), and $\big(\frac{p}{a},\frac{1-p}{b}\big)$ do not satisfy  item (i) or item (iii).

In fact, using that $C^1$ and $C^2$ are given by \eqref{stochmatr},
we have that  $0\leq e^{c^i_{1j}}\leq 1$, $a=e^{c^1_{11}}-e^{c^1_{12}}$ and $b=e^{c^2_{11}}-e^{c^2_{12}}$,  and this implies $-1\leq a\leq 1$ and $-1\leq b\leq 1$, and as a consequence $a-1\leq 0$ and $b-1\leq 0$. Then, if $a\neq b$, we have that $z_1=\frac{1-b}{a-b}$ and $z_2=\frac{a-1}{a-b}$ have opposite signs, and this  not satisfy item (i).

Let us now analyze the case $z_1=\frac{p}{a},z_2=\frac{1-p}{b}$. We will show that, either this solution does not satisfy item (i), or   the dominant eigenvalue of the matrix in \eqref{dominant} is greater than 1, and therefore do not satisfy item (iii). In fact, the matrix in \eqref{dominant} is
$$B=\frac{p}{e^{c^1_{11}}-e^{c^1_{12}}}\left(\begin{array}{cc}e^{c^1_{11}} & e^{c^1_{12}}\\ 1-e^{c^1_{11}} & 1- e^{c^1_{12}} \end{array}\right)+\frac{1-p}{e^{c^2_{11}}-e^{c^2_{12}}}\left(\begin{array}{cc}e^{c^2_{11}} & e^{c^2_{12}}\\ 1-e^{c^2_{11}} & 1- e^{c^2_{12}} \end{array}\right).$$
This matrix have eigenvalues given by $\lambda_1=1$ and $\lambda_2=\frac{p}{e^{c^1_{11}}-e^{c^1_{12}}}+\frac{1-p}{e^{c^2_{11}}-e^{c^2_{12}}}=\frac{p}{a}+\frac{1-p}{b}$. If this solution satisfy item (i), then $0<e^{-\varphi_1}=z_1=\frac{p}{a}$ and $0<e^{-\varphi_2}=z_2=\frac{1-p}{b}$, and then we have $0<a\leq 1$ and $0<b\leq 1$. This implies $\frac{p}{a}+\frac{1-p}{b}\geq p+1-p=1$, (with equality only if $a=b=1$, i.e., $z_1=p$ and $z_2=1-p$). In the case we have strict inequality, $\lambda_2>1$ will be the dominant eigenvalue.

We conclude that $(p,1-p)$ is the only critical point that satisfy items (i), (ii), and (iii). Hence, as we know that there exist a solution to 
$$\inf_{\varphi:P(c-\varphi)= 0} \int_{X} \varphi(x) \, d\mu $$ it must be given by 
$(\tilde{\varphi}_1,\tilde{\varphi}_2)=(-\log (p),
-\log (1-p))$, which finishes the proof of item b) of Theorem \ref{varphi-custoduascoordenadas}.
\cqd

\bigskip

\textbf{Remark:} Note that the original problem is to minimize $\int_{X} \varphi(x) \, d\mu =\varphi_1 p+\varphi_2(1-p)$ subject to the restriction
$\{\varphi:P(c-\varphi)= 0\} $. We know by \cite{LMMS2}, that exists a unique minimizer to this problem, that we denote by $\tilde{\varphi}=(\tilde{\varphi_1},\tilde{\varphi_2})$, and that  satisfies $P(c-\tilde{\varphi})= 0 $. Hence we get that the matrix $B_{\tilde{\varphi}}=e^{-\tilde{\varphi_1}}C^1+e^{-\tilde{\varphi_2}}C^2$ (see equation \eqref{dominant}) has 1 as dominant eigenvalue and this eigenvalue is simple. 

Instead of finding directly the minimizer of the original problem, we solved a second problem, which is to minimize $f(z_1,z_2)=- p \log z_1 -(1-p)\log z_2$, subject to $g(z_1,z_2)=0$, where $z_1=e^{-\varphi_1}, z_2=e^{-\varphi_2}$.

\medskip

\noindent\textbf{ Claim:} The solution of the original problem $\tilde z=(\tilde{z_1},\tilde{z_2})=(e^{-\tilde{\varphi_1}},e^{-\tilde{\varphi_2}})$ can be found in the set of extremal points of the second problem, that we determine using Lagrange multipliers.

\medskip

\noindent \textit{Proof of the claim:}

The conic $g(z)=0$ is given by points whose associated matrices have 1 as one of its eigenvalues. Such conic can have two subsets: The first  of then given by matrices having 1 as maximal eigenvalue (call $C_1$ this subset) and other (call it $C_2$) having 1 as minimal eigenvalue. ($C_1$ and $C_2$ can intersect in points whose associated matrices have both eigenvalues given by 1.)
Let us show that  $\tilde{z}$ is the interior of $C_1$. In fact, first note that  $B_{\tilde{z}}$ is stricly positive, which implies the dominant eigenvalue of  the matrix $B_{\tilde{z}}=\tilde{z_1}C^1+\tilde{z_2}C^2$ is simple, and  therefore the other eigenvalue is smaller than 1. 
Now, using the fact that the spectrum of $B_{{z}}={z_1}C^1+{z_2}C^2$ varies continuously in the parameter $z$, we know that there exists a neighbourhood $V_{\tilde{z}}$ of 
$\tilde{z}$ in the curve $g(z)=0$ such that the matrix $B_z$ has 1 as dominant and simple eigenvalue for any $z\in V_{\tilde{z}}$. This implies that  $\tilde{z}$ is the interior of $C_1$.
Finally, as $f(\tilde{z})<f(z)$ for all $z\in V_{\tilde{z}}$ (because $\tilde{z}$ is a minimum point of the original problem), this implies $\tilde{z}$ is a local minimum of the second problem and can be found using Lagrange multipliers.
{\it End of Proof of Claim.}

{\it Remark on the claim: }Note that the level curves of the function $f(z_1,z_2)$ have a simple geometry: $z_2$ is, in fact, proportional to a negative power of $z_1$, in any level curve. Also, $g(z_1,z_2)=0$ determines a conic, which also has a simple geometry, and this prevent any pathologies that could invalidate the argument above. 

\bigskip

\vspace{.5cm}
Now we want to perform Step 2.

The normalized cost associated to $A=c-\tilde{\varphi}$ is
$$\bar{A}=A+\log(h_A)-\log(h_A\circ \sigma)-\log(\lambda_A), $$
where $h_A$ is the eigenfunction  of $L_A=L_{c-\tilde{\varphi}}$ associated to the maximal eigenvalue $\lambda_A=\lambda_{c-\tilde{\varphi}}$. We see in the proof of the theorem \ref{varphi-custoduascoordenadas} that $\lambda_{c-\tilde{\varphi}}=1.$
And using \eqref{autofuncao} we get that $h_A=v_{c-\tilde{\varphi}}$ is the  left eigenvector of $B$ associated to the eigenvalue $1$, where 
$$B=\left(\begin{array}{cc}         e^{c^{1}_{11}} \,e^{-\tilde\varphi_1} +e^{c^{2}_{11}}\,e^{-\tilde\varphi_2} & e^{c^{1}_{12}} \,e^{-\tilde\varphi_1} +e^{c^{2}_{12}}\,e^{-\tilde\varphi_2}\\ e^{c^{1}_{21}}\,e^{-\tilde\varphi_1}+ e^{c^{2}_{21}}\,e^{-\tilde\varphi_2} & e^{c^{1}_{22}}\,e^{-\tilde\varphi_1}+ e^{c^{2}_{22}}\,e^{-\tilde\varphi_2} \end{array}\right)=z_1C^1+z_2C^2,
$$
i.e., $h_AB=h_A$, where $(z_1,z_2)=(e^{-\tilde{\varphi}_1},e^{-\tilde{\varphi}_2})$.

\bigskip

\textbf{Remark: }If we are in case b) of Theorem \ref{varphi-custoduascoordenadas}, then the cost $A=c-\tilde{\varphi}$ is already normalized.
In fact, if we denote by $p(1)=p=e^{-\tilde{\varphi}_1}$ and $p(2)=1-p=e^{-\tilde{\varphi}_2}$, using \eqref{stochmatr} we have
\[\sum_{x\in X} \sum_{\sigma(w)=y}e^{c(x,w)-\tilde{\varphi}(x)} =\sum_{x\in X}p(x) \sum_{a=1}^2e^{c(x,a,y_1)}=\sum_{x\in X}p(x)\Big[ e^{c(x,1,y_1)}+1-e^{c(x,1,y_1)}\Big]=1,\]
for all $y\in \Omega$.

\vspace{0.7cm}

\textbf{Example:} Let us suppose that the cost $c(x,y)=c(x,y_1,y_2)$ is represented by the following matrices
\[C^{1}=\left(\begin{array}{cc}3 & 5\\ 2 & 4 \end{array}\right)\,\,\,\,,\,\,\,\,C^{2}=\left(\begin{array}{cc}2 & 1\\ 4 & 3 \end{array}\right).\]

In order to perform Step 1 we look for positive solutions of \eqref{eqML}, and we get that the  unique positive solution  is $(z_1,z_2)=(0.101972,0.0568922)$.

Now in order to perform Step 2  we need first to calculate the left eigenvector, associated to the eigenvalue 1, of the matrix $B=z_1C^1+z_2C^2$.

We have that $$B=\left(\begin{array}{cc} 0.4197 & 0.566751\\ 0.431512 &0.578563 \end{array}\right)$$
and we can get that  $h_A=(0.596709,0.802458)$ is such that $h_A B=h_A$.

Then the  matrices 
$$\bar C^{1}=\left(\begin{array}{cc}0.3059 & 0.379132\\ 0.274264 & 0.407887 \end{array}\right)\,\,\,\,,\,\,\,\,\bar C^{2}=\left(\begin{array}{cc}0.113784 & 0.0423052\\ 0.306036 & 0.170677 \end{array}\right),$$
which  are obtain by the expression 
$$\bar{C}^x_{ij}=z_xC^x_{ij} \frac{h_A(i)}{h_A(j)}=e^{c(x,i,j)-\tilde{\varphi}(x)+\log(h_A(i))-\log(h_A(j))},$$
represent the normalized cost $\bar{A}=A+\log(h_A)-\log(h_A\circ \sigma)$.

In fact,  by the equation \eqref{defmatrix}, we have that 
$C^x_{ij}=e^{c(x,i,j)}$,  as $(z_1,z_2)=(e^{-\tilde{\varphi}(1)},e^{-\tilde{\varphi}(2)})$ we have $z_xC^x_{ij}=e^{c(x,i,j)-\tilde{\varphi}(x)}$, and finally  if $y=(i,j,y_3,y_4,....)$ then $\sigma(y)=(j,y_3,y_4,....)$. Therefore
$$\bar{C}^x_{ij}=e^{c(x,i,j)-\tilde{\varphi}(x)+\log(h_A(y))-\log(h_A\circ\sigma(y))}=e^{\bar{A}(x,y)}.$$

\bigskip
{\sc Joana Mohr\\
Departamento de Matematica\\
Universidade Federal do Rio Grande do Sul - UFRGS\\
Avenida Bento Gonçalves, 9500\\
91509-900 - Porto Alegre RS\\
Brazil}\\
{\it email joana.mohr@ufrgs.br}

\bigskip
{\sc Rafael Rigão Souza\\
Departamento de Matematica\\
Universidade Federal do Rio Grande do Sul - UFRGS\\
Avenida Bento Gonçalves, 9500\\
91509-900 - Porto Alegre RS\\
Brazil}\\
{\it email rafars@mat.ufrgs.br} 

\bigskip

{\footnotesize The author RRS was partially supported by FAPERGS (proc.002063-2551/13-0).}

\end{document}